\documentclass[12pt]{amsart}
 \usepackage[latin1]{inputenc}
 \usepackage[dvips]{graphicx}
\usepackage{hyperref}
 
 \usepackage{amsmath}
 \usepackage{amsthm}
 \usepackage{amsfonts}
 \usepackage{amssymb}
 \usepackage{amsopn}
 \usepackage{layout}
 \usepackage{verbatim}
 \usepackage{alltt}

 \input xypic
 \xyoption{all}

\newcommand{\del}[2]{{}}

\providecommand{\norm}[1]{\lVert#1\rVert}

\newcommand{\C}{\mathbb C}

\newcommand{\R}{\mathbb R}
\newcommand{\Q}{\mathbb Q}
\newcommand{\Z}{\mathbb Z}
\newcommand{\F}{\mathbb F}

\def\O{{\mathcal O}}
\def\a{{\mathfrak a}}
\def\q{{\mathfrak q}}
\def\p{{\mathfrak p}}

\def\H{{\mathbb H}}
\def\h{{\mathcal H}}
\def\sl{{\rm PSL_2}}   
\def\gl{{\rm PGL_2}}
\def\slo{{\rm PSL_1}}   
\def\glo{{\rm PGL_1}}

\title{Some Applications of Number Theory to 3-Manifold Theory}
\author{Mehmet Haluk \c{S}eng\"un}
\email{mehmet.sengun@uni-due.de}
\urladdr{http://www.uni-due.de/\char126hm0074/}
\address{Max Planck Institute for Mathematics, Bonn, Germany}

\begin{document}
 \maketitle
 \tableofcontents
\section{Introduction}
These are the extended notes of a talk I gave at the Geometric Topology Seminar of the Max Planck Institute for Mathematics in Bonn 
on January 30th, 2012. My goal was to familiarize the topologists with the basics of arithmetic hyperbolic 3-manifolds and sketch some interesting results in the theory of $3$-manifolds (such as Labesse-Schwermer \cite{labesse-schwermer}, Calegari-Dunfield \cite{calegari-dunfield}, Dunfield-Ramakrishnan \cite{dunfield-ramakrishnan}) that are obtained by exploiting the connections with number theory and automorphic forms. The overall intention 
was to stimulate interaction between the number theorists and the topologists present at the Institute. 

These notes are merely expository. I assumed that the audience had only a minimal background in number theory and I tried to 
give detailed citations referring the reader to sources. For many parts, I followed extensively 
the excellent books of Maclachlan-Reid \cite{maclachlan-reid} and Elstrodt-Grunewald-Mennicke \cite{egm}. Proofs of most of the facts 
on arithmetic lattices that I collected here can be found in these sources. I also benefited greatly 
from the expository articles of Reid \cite{reid-survey} and Schwermer \cite{schwermer}.  To make things more concrete and more 
accessible to the non-number theorist audience, I adopted the perspective of ``automorphic forms' as opposed to ``automorphic representations", 
which is a more suitable medium for the specific considerations in this paper and is the perspective adopted in most of the above cited literature. 

Please feel free to send me any comments.

\section{Quaternion Algebras}
The definition of arithmetic hyperbolic 3-manifold relies on the notion of quaternion algebra. In this section we collect 
results about quaternion algebras that we will need. 

A {\bf quaternion algebra} $D$ over a field $K$ (denoted $D / K$) is  a necessarily non-commutative ring with an injective ring 
map $K \rightarrow D$ such that 
\begin{enumerate}
\item the image of $K$ is the centre of $D$,
\item the dimension of $D$, considered as a $K$-vector space, is 4,
\item $D$ has no non-trivial two-sided ideals.
\end{enumerate}

It is well-known that if the characteristic of $K$ is not 2, then every quaternion algebra is isomorphic to one given in the following form:
$$\left ( \dfrac{a,b}{K}\right ) := \{ 1K \oplus iK \oplus jK \oplus ijK \mid i^2=a1, \ j^2=b1, \ ij=-ji \}. $$
Here the multiplication rules are extended via linearity to the whole vector space. 
The canonical example of a quaternion algebra is the $2\times 2$ matrix algebra $M_2(K) \simeq (\frac{1,1}{K})$. 
A well-known example is the so called {\it Hamiltonians} $\H \simeq (\frac{-1,-1}{\R})$.

In fact, every quaternion algebra over a field $K$ is either a skew-field (a.k.a. division ring) or isomorphic to $M_2(K)$. If $K$ is algebraically closed, then the only quaternion algebra $D / K$ is, up to isomorphism, $M_2(K)$. This implies that when we base change a quaternion algebra $D / K$ to the algebraic closure $\overline{K}$ of $K$, it will necessarily become isomorphic to $M_2(\overline{K})$. Hence, we may regard $D/K$ as a subalgebra of the matrix algebra $M_2(\overline{K})$.  In particular, we can define trace and norm maps on $D/K$ 
using the trace and determinant maps on $M_2(\overline{K})$. In fact, one does not need to go to the algebraic closure to 
get to the matrix algebra, it is known that there is always a degree two extension $L/K$ such that $D \otimes_K L \simeq M_2(L)$. 

When $K=\R$, there are, up to isomorphism, two quaternion algebras $D/K$, namely 
$\H$ and $M_2(\R)$. Similarly, when $K$ is a finite extension of $\Q_p$, there are again, up to isomorphism, only two, namely $M_2(K)$ and a skew-field which can be explicitly given (but we will not here) using the unique unramified degree 2 extension of $K$. 

Given a quaternion algebra $D$ over a number field $K$, we say that $D$ {\bf ramifies} over a place $v$ of $K$, if 
the base-change quaternion algebra $D \otimes_K K_v$, where $K_v$ is the completion of $K$ at the place $v$, is a 
skew-field. Note that if $v : K \hookrightarrow K_v$ is a place of $K$ and $D/K \simeq (\frac{a,b}{K})$, then 
$$D \otimes_K K_v \simeq \left ( \dfrac{v(a),v(b)}{K_v} \right ).$$ 
Observe that $D/K$ cannot ramify over a complex place $v$ since then $K_v \simeq \C$ and there is only 
one quaternion algebra over $\C$. On the contrary, ramification behaviour may vary among the real places. 
The {\bf discriminant} $\delta(D)$ of $D/K$ is the product of all prime ideals of $\O_K$ corresponding to the 
finite places of $K$ over which $D$ ramifies. Let us call the set of places of $K$ over which $D/K$ ramifies {\bf S(D)}. It turns out that $S(D)$ provides all the information one needs about $D/K$. It is well known that
\begin{enumerate}
\item S(D) has finite even cardinality,
\item if $S(D)=S(E)$, then $D/K \simeq E/K$,
\item for any set $S$ of places of $K$ of even cardinality that contains no complex place, there is a quaternion algebra 
$D/K$ such that $S(D)=S$. 
\end{enumerate}
Hence there is a bijection between quaternion algebras $D/K$ and finite sets $S$ of places of $K$ of even size which 
do not contain complex places.

An {\bf order} $\O$ in a quaternion algebra $D$ over a number field $K$ with ring of integers $\O_K$ is a finitely 
generated $\O_K$-module in $D$ such that $\O \otimes_{\O_K}K = D$ which is also a ring with 1. The number of conjugacy classes of 
maximal (with respect to inclusion) order of $D/K$ is finite. In particular, if $D/K$ is a skew-field, the elements of integral 
norm form a maximal order of $D$ and it is the unique maximal order up to conjugation.

\section{Hyperbolic Spaces}
For the purposes of this paper, we will need to consider the hyperbolic spaces of dimension 2 and 3. 
These arise as the symmetric spaces of the real Lie groups $\sl(\R)$ and $\sl(\C)$ respectively.

\subsection{Hyperbolic 2-space}
The hyperbolic 2-space $\h_2$ is the unique connected, simply connected Riemannian manifold of dimension 2 with 
constant sectional curvature $-1$.  A standard model for $\h_2$ is the upper half plane model 
$$\{ (x,y) \in \R \times \R \mid y>0 \}$$
with the metric coming from the line element
$$ds^2 = \dfrac{dx^2+dy^2}{y^2}$$
under which the distance $d(P,P')$ between two points $P=(x,y)$, $P'=(x',y')$ is given by
$$\cosh d(P,P')= 1+\dfrac{(x-x')^2+(y-y')^2}{2yy'}.$$

The orientation-preserving isometries of $\H$ can be identified with the group $\sl(\R)$. The action on $\h_2$ is given as 
$$\begin{pmatrix} a & b \\ c&d \end{pmatrix} \cdot z =\dfrac{az+b}{cz+d}$$
for $z=x+iy \in \h_2$. Note that 
$${\rm Re} \left (\dfrac{az+b}{cz+d} \right ) = \dfrac{(ax+b)(cx+d)+acy^2}{(cx+d)^2+c^2y^2}, \ \ \ \ {\rm Im} \left (\dfrac{az+b}{cz+d} \right)=\dfrac{y}{(cx+d)^2+c^2y^2}.$$
 
\subsection{Hyperbolic 3-space}
Hyperbolic 3-space $\h_3$ is the unique connected, simply connected Riemannian manifold of dimension 3 with 
constant sectional curvature $-1$.  A standard model for $\h_3$ is the upper half space model 
$$\{ (x,y) \in \C \times \R \mid y>0 \} = \{ (x_1,x_2,y) \in \R^3 \mid y>0 \}  $$
with the metric coming from the line element
$$ds^2 = \dfrac{dx_1^2+dx_2^2+dy^2}{y^2}$$
under which the distance $d(P,P')$ between two points $P=(x_1,x_2,y)$,$P'=(x'_1,x'_2,y')$ is given by
$$\cosh d(P,P')=1+ \dfrac{(x_1-x'_1)^2+(x_2-x'_2)^2+(y-y')^2}{2yy'}.$$

The orientation-preserving isometries of $\h_3$ can be identified with the group $\sl(\C) \simeq \gl(\C)$. The action on $\h_3$ is given as 
$$\begin{pmatrix} a & b \\ c&d \end{pmatrix} \cdot z = 
\Biggl ( \frac{(ax+b)\overline{(cx+d)}+a\bar{c}y^2}{|cx+d|^2+|c|^2y^2},\frac{y}{|cx+d|^2+|c|^2y^2} \Biggr )$$
where $z=(x,y) \in \h_3$. 

It is computationally convenient to regard $\h_3$ inside the Hamiltonians $\H$. 
We embed $\h_3$ into $\H$ via $(x_1,x_2,y) \mapsto (x_1,x_2,y,0).$ In other words, $(x,y) \mapsto x+y j$. 
In particular, the element $j \in \H$ corresponds to $(0,0,1) \in \R^3$. Then the action of 
$\sl(\C)$ on $\h_3$ takes the familiar form
$$\begin{pmatrix} a & b \\ c&d \end{pmatrix} \cdot z = (az+b)(cz+d)^{-1} $$
where the inverse of $cz+d$ is taken inside the skew-field $\H$. 

To see this transition, consider $\H$  as a subalgebra of $\H \otimes_\R \C \simeq M_2(\C)$ as follows
$$\H \hookrightarrow M_2(\C), \ \ \ x+yj \mapsto \begin{pmatrix} x & -y \\ \bar{y} & \bar{x}\end{pmatrix}.$$
In particular,  if $x \in \C$, we have $x \mapsto (\begin{smallmatrix} x & 0 \\ 0 & \bar{x}\end{smallmatrix})$. This embedding respects the ring operations of $\H$. It is now trivial to verify the equivalence between the two descriptions of the action of $\sl(\C)$ on $\h_3$ that we gave above when we consider this matrix representation of $\h_3 \subset \H$.
 
\section{Arithmetic Lattices}

Consider a quaternion algebra $D/K$ over a number field $K$ of degree $n$ with $s$ complex places. Let $r \geq 0$ be the number of real places of $K$ that $D$ does not ramify over. We will need the groups $\glo(D)$ and $\slo(D)$, the group of invertibles elements of $D$ modulo its center and the goup of norm one elements modulo its center, respectively. For example, we have $\slo(M_2(K))\simeq \sl(K)$ and $\slo(\H) \simeq SU(2) / \{ \pm 1\}  \simeq SO(3)$, a compact Lie group.  Let $\O$ be an order in $D$. Then $\slo(\O)$ is discrete when considered as a subgroup
$$ \slo(\O) \hookrightarrow \prod_{v \ {\rm infinite}} \slo(D \otimes_K K_v) = \sl(\C)^s \times \sl(\R)^r \times SO(3)^{n-2s-r}.$$  
As the other factors are compact, it can be shown that the image of $\slo(\O)$ under the projection map onto the factor 
$\sl(\C)^s \times \sl(\R)^r$ is a discrete subgroup of $\sl(\C)^s \times \sl(\R)^r$. The map from $\slo(\O)$ onto its discrete 
image in $\sl(\C)^s \times \sl(\R)^r$ has finite kernel. We will, for convenience, denote this image with $\slo(\O)$ as well.

When $K$ is totally real and $D/K$ ramifies at all real places but one, $\slo(\O)$ gives a discrete subgroup of $\sl(\R)$ 
for every order $\O \subset D$. Especially important for us is the case when $K$ has exactly one complex place and $D/K$ 
ramifies over all the real places of $K$. Then $\slo(\O)$ gives a discrete subgroup of $\sl(\C)$ for every order $\O \subset D$. 

A discrete subgroup $\Gamma$ of $ \sl(\C)^s \times \sl(\R)^r $ is called a {\em lattice} if the associated orbifold (manifold is $\Gamma$ is 
torsion-free) $\Gamma \backslash \h_3^s \times \h_2^r$ has finite volume (we say that $\Gamma$ has finite covolume). Two discrete subgroups $\Gamma_1$,$\Gamma_2$ are {\em commensurable} (denoted $\Gamma_1 \equiv \Gamma_2$) if their intersection $\Gamma_1 \cap \Gamma_2$ is finite index in both $\Gamma_1$ and $\Gamma_2$. We say that they are {\em widely 
commensurable} if they are commensurable after possibly conjugation. We define the {\bf commensurator} $Comm(\Gamma)$ of a discrete subgroup $\Gamma$ as 
$$Comm(\Gamma) = \{ \gamma \in  \sl(\C)^s \times \sl(\R)^r  \mid \gamma\Gamma \gamma^{-1} \equiv \Gamma \}. $$
This is a group that clearly contains $\Gamma$. 

A discrete subgroup $\Gamma$ of $ \sl(\C)^s \times \sl(\R)^r $ is called {\bf arithmetic} if it is widely commensurable with some discrete subgroup $\slo(\O)$ that arises in the above way from some order $\O$ in some quaternion algebra $D/K$ over some number field $K$. In this case, we call $K$ the {\em field of definition} of $\Gamma$ and call $D/K$ the {\em defining quaternion algebra}.  It follows from a general result in the theory of arithmetic groups that every arithmetic subgroup is a lattice, that is, has finite covolume. Moreover, the associated orbifold is non-compact if and only if $S(D) =\emptyset$.

In the above definition, if we take $K$ totally real and $D/K$ ramified at all real places but one, we get {\bf arithmetic Fuchsian groups}. In this case, $\Gamma \subset \sl(\R)$ is non-cocompact if and only if $K=\Q$ and $D/K \simeq M_2(\Q)$. The case 
when $K$ has only one complex place and $D/K$ ramifies at all real places leads to {\bf arithmetic Kleinian groups}. In this case, 
$\Gamma \subset \sl(\C)$ is non-cocompact if and only if $K$ is imaginary quadratic and $D/K \simeq M_2(K)$.

The simplest examples of arithmetic Kleinian groups are the so called {\em Bianchi groups}. Take $K$ to be imaginary quadratic and $D/K \simeq M_2(K)$. Then $M_2(\O_K)$ is a maximal order and $\slo(M_2(\O_K)) \simeq \sl(\O_K)$ is the Bianchi group associated 
to $K$. As mentioned above, any non-cocompact arithmetic lattice in $\sl(\C)$ is widely commensurable with a Bianchi group.

Another well-known example of an arithmetic lattice is the complement of the figure eight knot. In fact, Alan Reid proved in \cite{reid} that this is the only arithmetic {\it knot complement}. Its fundamental group is isomorphic to an index 12 subgroup of the Bianchi group $\sl(\O_K)$ associated to $K=\Q(\sqrt{-3})$. 

Other well-known examples of arithmetic lattices are provided by the {\it hyperbolic tetrahedral groups}. A hyperbolic tetrahedral group is the index 2 subgroup consisting of orientation-preserving isometries in the discrete group generated by reflections in the faces of a hyperbolic tetrahedron (with possible ideal vertices) whose dihedral angles are submultiples of $\pi$. Lann\'er proved in 1950 \cite{lanner} that there are 32 such hyperbolic tetrahedra. It is known that 23 of these groups are arithmetic. They are discussed in detail in \cite{maclachlan-reid, egm}.

\section{Automorphic Forms}
\subsection{Forms on $\mathcal{H}_2$}
Let us introduce the automorphic forms that are associated to arithmetic lattices in $\sl(\R)$. 
These are complex valued real analytic functions on $\H$ with certain transformation properties which 
satisfy certain differential equations and growth properties.

Given $\gamma=(\begin{smallmatrix} a & b \\ c&d \end{smallmatrix}) \in \sl(\C)$ and $z \in \h_2$, 
let us introduce the multiplier system
$$J(\gamma, z):= (cz+d)$$

Given a function $F: \h_2 \rightarrow \C$ and $\gamma \in \sl(\R)$, we define the {\em slash operator} 
$$(F |_k\gamma)(z):=J(\gamma, z)^{-k} F(\gamma z)$$
for $z \in \h_2$.

The center of the universal enveloping algebra of the Lie algebra associated to the real Lie group $\sl(\R)$ is 
generated by an element $\Psi$ which acts on real analytic functions $F: \h_2 \rightarrow \C$ as a differential 
operator. 

Let $\Gamma$ be an arithmetic lattice in $\sl(\R)$ with defining field $K$ and defining quaternion algebra $D/K$. An {\bf automorphic form} for $\Gamma$ with weight $k$ and eigenvalue $\lambda$ is a real analytic function $F: \h_2 \rightarrow \C$ with the following properties.
\begin{enumerate}
\item $F|_k\gamma = F$ for every $\gamma \in \Gamma$,
\item $\Psi F = \lambda F$,
\item if $\Gamma$ is non-cocompact, then $F$ has at worst polynomial growth at each cusp.
\end{enumerate}

The set $M(\Gamma,k,\lambda)$ of automorphic forms for $\Gamma$ with weight $k$ and eigenvalue 
$\lambda$ is a finite dimensional complex vector space. 

A {\bf weight $2$ cuspidal modular form} for $\Gamma$ is a holomorphic function $F : \h_2 \rightarrow \C$ 
with the following properties.
\begin{enumerate}
\item $F(z)dz$ is a holomorphic differential 1-form on $\h_2$ that is $\Gamma$-invariant,
\item If $\Gamma$ is non-cocompact, then $\int_{\R / \Z} (F|_2\gamma )(x,y) dx = 0$ for every $\gamma \in \sl(\Z).$
\end{enumerate}
The last condition, in which case $K$ is necessarily $\Q$, is equivalent to saying that the constant coefficient in the {\em Fourier expansion} of $F|\gamma$ is equal to zero for every $\gamma \in \sl(\Z)$. Let us explain this. As $\Gamma$ is non-cocompact, it contains parabolic elements. The $\Gamma$-invariance of $F$, which is implicit in condition (1), implies that $F$ is a periodic function. It follows that the $F$ has a Fourier expansion of the form 
$$F(z)=\sum_{n>0} a_n e^{2\pi i nz/h}.$$

The space of weight 2 cuspidal modular forms for a fixed $\Gamma$ form a finite dimensional complex vector space which we will denote with $S_2(\Gamma)$. As automorphic forms for $\Gamma$, we have $S_2(\Gamma) \subset M(\Gamma,2,0)$.

\subsection{Forms on $\mathcal{H}_3$}
Let us introduce the automorphic forms that are associated to arithmetic lattices in $\sl(\C)$. 
These are {\em vector valued} real analytic functions on $\h_3$ with certain transformation properties satisfying 
certain differential equations and growth properties. The fact that our automorphic forms on $\h_2$ were 
scalar valued whereas automorphic forms on $\h_3$ are vector valued comes from the following fact: 
$SO(2)$, the maximal compact subgroup of $\textrm{SL}_2(\R)$, is abelian and hence has only one dimensional 
irreducible complex representations, on the other hand, $SU(2)$, the maximal compact subgroup of $\textrm{SL}_2(\C)$, 
is non-abelian and it has irreducible complex representations of every degree.

Given $\gamma=(\begin{smallmatrix} a & b \\ c&d \end{smallmatrix}) \in \sl(\C)$ and $z=x+yj \in \h_3$, 
let us introduce the multiplier system
$$J(\gamma, z):= \begin{pmatrix} cx+d & -cy \\ \bar{c}y & \overline{cx+d}\end{pmatrix}$$

Given a function $F: \h_3 \rightarrow \C^{k+1}$ and $\gamma \in \sl(\C)$, we define the {\em slash operator} 
$$(F |_k\gamma)(z):=\sigma^k(J(\gamma, z)^{-1}) F(\gamma z)$$
where $\sigma^k$ is the {\em symmetric $k^{th}$ power} of the standard representation of $\sl(\C)$ on $\C^2$.

The case $k=2$ will be especially important for us. In this case we have $F: \h_3 \rightarrow \C^3$ and
$$(F |_k\gamma)(z)=\dfrac{1}{|r|^2+|s|^2} 
\begin{pmatrix}  \bar{r}^2 & 2\bar{r}s & s^2 \\ 
                          -\bar{r}\bar{s} & |r|^2-|s|^2 & rs \\ 
                          \bar{s}^2 & -2r\bar{s} & r^2 \end{pmatrix} F(\gamma z)$$
where $\gamma=( \begin{smallmatrix} a &b \\ c&d \end{smallmatrix})$ and $r=cx+d$ and $s=cy$.

The center of the universal enveloping algebra of the Lie algebra associated to the real Lie group $\sl(\C)$ is 
generated by two elements $\Psi, \Psi'$. These act on real analytic functions $F: \h_3 \rightarrow \C^{k+1}$ as differential 
operators. 

Let $\Gamma$ be an arithmetic lattice in $\sl(\C)$ with defining field $K$ and defining quaternion algebra $D/K$. An {\bf automorphic form} for $\Gamma$ with weight $k$ and eigenvalues $(\lambda, \lambda')$ is a real analytic function $F: \h_3 \rightarrow \C^{k+1}$ with the following properties.
\begin{enumerate}
\item $F|_k\gamma = F$ for every $\gamma \in \Gamma$,
\item $\Psi F = \lambda F$ and $\Psi' F = \lambda' F$,
\item if $\Gamma$ is non-cocompact, then $F$ has at worst polynomial growth at each cusp.
\end{enumerate}

The set $M(\Gamma,k,\lambda,\lambda')$ of automorphic forms for $\Gamma$ with weight $k$ and eigenvalues $(\lambda, \lambda')$ is a finite dimensional complex vector space. 

As there is no complex structure on $\h_3$, we will replace holomorphicity with harmonicity. Let $\beta_1:=-\frac{dx}{y}, \beta_2:= \frac{dy}{y}, \beta_3:=\frac{d\bar{x}}{y} $ be a basis of differential 1-forms on $\h_3$. 
A differential form $\omega$ is {\em harmonic} if $\Delta \omega =0$ where $\Delta=d \circ \delta + \delta \circ d$ 
is the usual Laplacian with $d$ being the exterior derivative and $\delta$ the codifferential operator. 
Then $\sl(\C)$ acts on the space of differential 1-forms as 
$$\gamma \cdot {}^t(\beta_1,\beta_2,\beta_3)_{(z)} = \sigma^2(J(\gamma,z)){}^t(\beta_1,\beta_2,\beta_3)_{(z)} .$$
A {\bf weight $2$ cuspidal modular form} for $\Gamma$ is a real analytic function $F=(F_1,F_2,F_3) : \h_3 \rightarrow \C^3$ 
with the following properties.
\begin{enumerate}
\item $F_1\beta_1 + F_2 \beta_2+F_3\beta_3$ is a harmonic differential 1-form on $\h_3$ that is $\Gamma$-invariant,
\item If $\Gamma$ is non-cocompact, then $\int_{\C / \O_K} (F|_2\gamma )(x,y) dx = 0$ for every $\gamma \in \sl(\O_K).$
\end{enumerate}
The last condition, in which case $K$ is necessarily an imaginary quadratic field, is equivalent to saying that the constant coefficient in the {\em Fourier-Bessel expansion} of $F|\gamma$ is equal to zero for every $\gamma \in \sl(\O_K)$. Let us explain this. As $\Gamma$ is non-cocompact, it contains parabolic elements. The $\Gamma$-invariance of $F$, which is implicit in condition (1), implies that $F$ is a periodic function in the $x=(x_1,x_2)$-variable. It follows that the $F$ has a Fourier-Bessel expansion of the form 
$$F(x,y)=\sum_{0 \not =\alpha \in \O_K}c(\alpha) y^2 \mathbb{K}\left ( \dfrac{4\pi|\alpha|y}{\sqrt{|\triangle|}} \right )
                        \psi\left (\dfrac{\alpha x}{\sqrt{\triangle}} \right )$$
where 
$$\psi(x)=e^{2\pi(x+\bar{x})}$$
and  
$$\mathbb{K}(t)=\left ( -\dfrac{i}{2}K_1(y),K_0(y),\dfrac{i}{2}K_1(y) \right)$$
with $K_0,K_1$ are the hyperbolic Bessel functions satisfying the differential equation 
$$\dfrac{dK_j}{dy^2}+\dfrac{1}{y}\dfrac{dK_j}{dy}-\left ( 1+\dfrac{1}{y^{2j}}\right )K_j = 0, \ \ \ \ j=0,1$$
and decreases rapidly at infinity. 

The space of weight 2 cuspidal modular forms for a fixed $\Gamma$ form a finite dimensional complex vector space which we will denote with $S_2(\Gamma)$. As automorphic forms for $\Gamma$, we have $S_2(\Gamma) \subset M(\Gamma,2,0,0)$.

\subsection{Forms on $\h_3^s \times \h_2^r$} 
For an arithmetic lattice $\Gamma \subset \sl(\C)^s \times \sl(\R)^r$, there is a {\em general notion} of weight 2 cuspidal modular form for $\Gamma$. 
As you can guess, these are vector valued functions on $\h_3^s \times \h_2^r$ which are ``weight 2 cuspidal modular forms" in every component. As we have explained what it means to be a weight 2 cuspidal modular form for both component types, we will not discuss this general notion any further. The interested reader can find very detailed general accounts by Hida \cite{hida, hida-2}. See also Friedberg \cite{friedberg1, friedberg2} and Bygott \cite{bygott} 
for a detailed study of automorphic forms on $\H_3$. We shall use same symbol $S_2(\Gamma)$ to denote the finite dimensional complex vector space of weight 2 cuspidal modular forms for such a lattice $\Gamma$.

These spaces $S_2(\Gamma)$ come equipped with an infinite family of commuting endomorphisms $\mathbb{T}$, called {\em Hecke operators}, which are indexed by elements in $\O_K$, the ring of integers of the field of definition of $\Gamma$. We will not define these operators here, instead we will define their counterparts in the cohomological setting. 

The existence of this family of operators is crucial for connections with number theory. To illustrate, assume that $f \in S_2(\Gamma)$ is a simultaneous 
eigenvector for the action of $\mathbb{T}$, that is, for each $T_{\alpha} \in \mathbb{T}$, $\alpha \in \O_K$, we have $T_{\alpha} f = a_f(\alpha) \cdot f$ for some $a_f(\alpha) \in \C$. It can be shown that the $a_f(\alpha)$ all live in a number field $L$ and in fact, when $f$ is suitably normalized, they live 
in $\O_L$, the ring of integers of $L$. Langlands Programme conjectures that there is an irreducible continuous representation 
of the absolute Galois group $G_K$ of the field of definition $K$ of $\Gamma$
$$\rho : G_K \rightarrow \textrm{GL}_2(\overline{\Q}_p)$$ 
with associated {\bf arithmetic data} $\{ a_{\rho}(\alpha) \}_{\alpha \in \O_K}$ such that $a_{\rho}(\alpha) = a_f(\alpha)$ for almost all $\alpha \in \O_K$. Here the ``arithmetic data" associated to $\rho$ is the result of an explicit procedure that we will not discuss.
 
\subsection{Jacquet-Langlands Correspondence}
The so called Jacquet-Langlands correspondence, see \cite{jacquet-langlands}, is a statement about automorphic representations of the algebraic group $\textrm{Res}_{K/\Q} \gl$ and its inner forms. We will present this correspondence in the setting of automorphic forms.

Let $K$ be a number field  with ring of integers $\O_K$. Let $D/K$ be a quaternion algebra with $S(D) \not= \emptyset$. Let $\O$ be some order in $D$ and $\Gamma$ be the arithmetic cocompact lattice given by $\slo(\O)$. Then there is an injection
$$S_2(\Gamma) \hookrightarrow S_2(\Gamma')$$
where $\Gamma'$ is a finite index subgroup of the $\sl(\O_K)$. Moreover, the above linear map respects the action 
of the Hecke operators on both spaces. In fact, the (two-way) passage from $\Gamma$ to $\Gamma'$ is completely explicit 
and the image of the linear map can be made described precisely using the theory of ``old/new"-forms. 
We will denote the image of this map with $S_2(\Gamma')^{\delta(D)\textrm{-new}}.$

\section{Hecke Operators on the Cohomology}

In this section we will consider the cohomology  groups $H^1(\Gamma \backslash \h_3, \Z)$ for lattices $\Gamma$ in $\sl(\C)$, 
even though all our discussion generalizes directly to the setting of general lattices in $\sl(\C)^s \times \sl(\R)^r$. 
When the lattice is arithmetic, there is a very special infinite family of operators, called Hecke operators, that act on these 
cohomology groups. We will study these operators.

Let $\Gamma$ be any lattice in $\sl(\C)$. Given $g \in Comm(\Gamma)$, let us consider the 3-folds 
$M,M_g,M^g$ associated to the lattices $\Gamma, \Gamma_g := \Gamma \cap g \Gamma g^{-1}, 
\Gamma^g :=\Gamma \cap g^{-1} \Gamma g $ 
respectively. We have {\em finite} coverings, induced by inclusion of fundamental groups, 
$$r_g:M_g \rightarrow M, \ \ \ \ \ r^g: M^g \rightarrow M$$
and an isometry 
$$\tau : M_g \rightarrow M^g$$
induced by conjugation by $g$ isomorphism between $\Gamma_g$ and $\Gamma^g$.
The composition $s_g := r^g \circ \tau$ gives us a second finite covering from $M_g$ to $M$. 
The coverings $r_g$ induce linear maps between the homology groups 
$$r_g^* : H^1(M,\Z) \rightarrow H^1(M_g,\Z).$$
The process of summing the finitely many preimages in $M_g$ of a point of $M$ under $s_g$ leads to 
$$s_g^* : H^1(M_g,\Z) \rightarrow H^1(M,\Z).$$
Note that $s_g^*$ is equivalent to the composition
$$H^1(M_g,\Z) \rightarrow H^1(M^g,\Z) \rightarrow H^1(M,\Z)$$
where the first arrow is induced by $\tau$, 
and the second arrow is simply the corestriction map (which corresponds to the transfer map of group cohomology). We define the {\bf Hecke operator} $T_g$ associated to $g \in Comm(\Gamma)$ as the composition 
$$T_g:= s_g^* \circ r_g^* : H^1(M,\Z) \rightarrow H^1(M,\Z).$$  
There is a notion of isomorphism of Hecke operators that we shall not present. It turns out that up to isomorphism, 
$T_g$ depends only on the double coset $\Gamma g \Gamma$.

One can define Hecke operators using the above process for cohomology and homology groups of every degree. 
It might provide insight to look at the situation from the perspective of $H_2(M,\Z)$, which is isomorphic to $H^1(M,\Z)$ 
by Poincar\'e duality when $M$ is compact. If an emdedded surface $S \subset M$ represents $\omega \in H_2(M,\Z)$, then the immersed 
surface $r_g(s_g^{-1}(S))$ represents $T_g\omega \in H_2(M,\Z)$. 

It is a classical result of Margulis that either $\Gamma$ is finite index in its commensurator or its commensurator 
is dense in $\sl(\C)$. Moreover, the latter happens if and only if $\Gamma$ is arithmetic.
The commensurator of a Bianchi group $\sl(\O_K)$ is $\gl(K)$. 
More generally, the commensurator of an arithmetic lattice $\Gamma$ in $\sl(\C)$ with defining quaternion 
algebra $D/K$ is $\glo(D) \subset \sl(\C)$.

Especially important is the subfamily $\mathbb{T}$ of Hecke operators associated to the prime elements of 
the ring of integers $\O_K$ of the definition field $K$. More precisely,  let $\pi$ be a prime element  of $\O_K$. 
Then we have $\left ( \begin{smallmatrix} \pi & 0 \\ 0 & 1 \end{smallmatrix} \right ) \in Comm(\Gamma)$ (note that 
here we regard $D/K \subset M_2(\overline{K})$). If $T_{\pi}$ is the Hecke operator associated to 
$\left ( \begin{smallmatrix} \pi & 0 \\ 0 & 1 \end{smallmatrix} \right )$, then we set 
$$\mathbb{T}=\{ T_{\pi} \mid \pi \in \O_K \}.$$

\section{Cohomology and Automorphic Forms}
Let $\Gamma$ be an arithmetic lattice in $\sl(\C)$ with associated $3$-fold $Y$. 

Assume that $Y$ is non-compact. Without loss of generality, we assume that $\Gamma$ is finite index subgroup of a Bianchi 
group. Borel and Serre proved in \cite{borel-serre} that there is a compact 3-fold $X$ with boundary (known as the {\em Borel-Serre compactification}), such that the interior of $X$ is homeomorphic to $Y$ and the embedding $Y \hookrightarrow X$ is a homotopy equivalence. This implies that $H^i(X,\C) \simeq H^i(Y,\C)$ for every $i \geq 0$. The boundary $\partial X$ of $X$ is a disjoint union of $2$-tori, each one corresponding to a cusp of $Y$. In fact, topologically, we can think of $X$ as the manifold obtained by attaching a $2$-torus at infinity to each cusp of $Y$. It is a classical result that the number of cusps of $Y$ when $\Gamma$ is the Bianchi group associated to $K$ is equal to the class number of $K$. 

Consider the map $res: H^1(X,\C) \rightarrow H^1(\partial X, \C)$ given by the restriction to the boundary.  The kernel of this map gives a subspace of $H^1(Y,\C)$ which is called the {\em cuspidal cohomology}, denoted $H^1_{cusp}(Y,\C)$. Harder showed in \cite{harder} that there is a section of the above restriction map which gives a subspace of $H^1(Y,\C)$ which is called the {\em Eisenstein cohomology} such that the decomposition 
$$H^1(Y,\C) = H^1_{cusp}(Y,\C) \oplus H^1_{Eis}(Y,\C)$$
is invariant under the action of the special family $\mathbb{T}$ of Hecke operators. Observe that when $Y$ is compact, we do not 
have the Eisenstein cohomology anymore, that is $H^1(Y,\C)=H^1_{cusp}(Y,\C)$. In a natural sense, while the Eisenstein cohomology comes from the contribution of the boundary of the Borel-Serre compactification $X$, the cuspidal cohomology belongs to the interior of $X$.  

While one has a good control over the Eisenstein cohomology, the cuspidal part is very mysterious. The importance of cuspidal cohomology comes from the fact that it can be identified with certain types of automorphic forms called ``cohomological". 

Now let us get back to the general case where $Y$ is not necessarily non-compact. Let $S_2(\Gamma)$ denote the space of weight $2$ cuspidal modular forms for an arithmetic lattice $\Gamma$ as discussed above. Then there is an isomorphism, called {\bf generalized Eichler-Shimura Isomorphism}
$$S_2(\Gamma) \simeq H^1_{cusp}(\Gamma,\C)$$
which respects the Hecke action on $S_2(\Gamma)$ and the action of $\mathbb{T}$ on $H^1_{cusp}(\Gamma,\C)$. The fact 
that this is not just an isomorphism of vector spaces but of Hecke modules is crucial for number theory. 

Roughly speaking, the above isomorphism comes from the facts that weight 2 cuspidal modular forms are
 essentially harmonic differential 1-forms and every cohomology class in the de Rahm cohomology 
$$H^1_{dR}(\Gamma \backslash \H, \R) \hookrightarrow H^1(\Gamma \backslash \H, \R) \simeq H^1(\Gamma,\R)$$
can be represented by a harmonic differential 1-form. The general connection between (vector valued) harmonic differential k-forms and cohomological automorphic forms was studied by Matsushima and Murakami, see \cite{matsushima-murakami}, in the cocompact setting.

All of the above was studied more generally by Harder in \cite{harder, harder2} for the algebraic group $\textrm{Res}_{K / \Q}(\gl)$ for any number field $K$. Let $\Gamma$ be an arithmetic lattice in $\sl(\C)^s \times \sl(\R)^r$ with associated quotient manifold $Y$. If $Y$ is non-compact, then we define the cuspidal and the Eisenstein parts of $H^{r+1}(Y,\C)$ using the Borel-Serre 
compactification of $Y$ as above. Again, if $Y$ is compact, the whole cohomology is cuspidal. The generalized Eichler-Shimura isomorphism 
tells us that 
$$S_2(\Gamma) \simeq H^{s+r}_{cusp}(Y,\C)$$
as Hecke modules.

The above isomorphism is the generalization of the classical Eichler-Shimura isomorphism between classical modular forms and the cohomology of arithmetic lattices in $\sl(\R)$. We should note that in this classical case, the cuspidal first cohomology of the lattice $\Gamma$ captures 
two copies of $S_2(\Gamma)$. In \cite{franke} J.Franke generalized this sort of connection to the fullest, that is, to the case of arithmetic lattices in general real Lie groups and their associated modular forms. Today most of the popular methods for computing with modular forms, such as the modular symbols method, is based on this passage to the (co)homology.

\subsection{Jacquet-Langlands Revisited: The Mystery}
When transfered to the setting of cohomology, our statement that was an unprecise version of the Jacquet-Langlands correspondence takes the following form. 

Let $K$ be a number field with $s$ complex places and $r$ real places. Let $D/K$ be a quaternion algebra ramifying over $0 \leq r_1$ of the real places of $K$. Let $\O$ be some order in $D$ and $\Gamma$ be the arithmetic cocompact lattice in $\sl(\C)^s \times \sl(\R)^{r-r_1}$ given by $\O$. 
Then there is an isomorphism
$$H^{s+r-r_1}(\Gamma,\C) \simeq H^{s+r}_{cusp}(\Gamma',\C)^{\delta(D)\textrm{-new}}$$
where $\Gamma'$ is a finite index subgroup of $\sl(\O_K)$. Moreover, the above linear map respects the action 
of the Hecke operator family $\mathbb{T}$ on both sides. 

As we mentioned before, the above relationship is born from trace formula comparisons in the realm of representation theory of adelic algebraic groups. From the perspective of topology, it is a {\bf complete mystery!} To elaborate on this mystery, let us consider the following two simplest cases.

Let $\Gamma \subset \sl(\R)$ be a cocompact arithmetic Fuchsian group with defining field $\Q$. Then we have 
$$H^1(\Gamma,\C) \hookrightarrow H^1(\Gamma',\C)$$
for some finite index subgroup $\Gamma' \subset \sl(\Z)$. 

Let $\Gamma \subset \sl(\C)$ be a cocompact arithmetic Kleinian group with an imaginary quadratic defining field $K$. Then we have 
$$H^1(\Gamma,\C) \hookrightarrow H^1(\Gamma',\C)$$
for some finite index subgroup $\Gamma' \subset \sl(\O_K)$. 

From a topological point of view, it is {\bf completely unknown} why the positiveness of the betti number of the compact manifold associated 
to $\Gamma$ implies the positiveness of the betti number of the non-compact manifold associated to $\Gamma'$.

\section{Volumes of Arithmetic Hyperbolic 3-Manifolds} 
By the Mostow-Prasad rigidity, the volume of a hyperbolic 3-manifold is a topological invariant. One of the attractive features of 
arithmetic hyperbolic 3-folds is that their volumes can be computed explicitly from the arithmetic data associated 
to the their fields of definition.

Let $K$ be a number field with ring of integers $\O_K$. Given an ideal $\p$ of $\O_K$, let $N\p$ denote its norm, which is the cardinality of the finite ring $\O / \p$. The {\bf Dedekind zeta function} associated to $K$ is defined as 
$$\zeta_K (s) :=\sum_{\p \triangleleft \O_K} \dfrac{1}{(N\p)^s}$$
which converges absolutely for complex numbers $s$ with $Re(s)>1$. 

Now assume that $K$ has exactly one complex place and $r \geq 0$ real places. Let $D/K$ be a quaternion algebra over $K$ ramifying over all real places. Let $F(D)$ denote the subset of $S(D)$ which consists of finite places of $K$. For each $v \in F(D)$, let $Nv$ denote the norm of the ideal $\p_v$ corresponding to the place $v$. Take a {\em maximal} 
order $\O$ in $D/K$ and consider the arithmetic lattice $\Gamma:=\slo(\O) \subset \sl(\C)$ as discussed above. 
Then the {\bf volume} $vol(\Gamma)$ of the hyperbolic 3-fold $\Gamma \backslash \H$ is 
$$vol(\Gamma) = \left ( \prod_{v \in F(D)} (Nv-1) \right ) \cdot \dfrac{|\bigtriangleup_K|^{3/2}}{(2\pi)^{2r+2}}\cdot \zeta_K(2)$$
where $\bigtriangleup_K$ is the discriminant of $K/\Q$.  In particular, we see that for an imaginary quadratic 
field $K$, the volume of the associated Bianchi group is given as 
$$vol(\sl(\O_K))=\dfrac{|\bigtriangleup_K|^{3/2}}{4\pi^2}\cdot \zeta_K(2)$$
as first proved by Humbert in 1919 in \cite{humbert}. A gap in his proof was filled by Grunewald and K\"uhnlein in 1998 in \cite{grunewald-kuhnlein}. For amusement, let us list the volumes of the 3-folds associated to some Bianchi groups. Let $\O_d$ be the ring of integers 
of the imaginary quadratic field $\Q(\sqrt{-d})$. Then 
\begin{center}
\begin{tabular}{c}
$vol(\sl(\O_1))\simeq0.305321$ \\
$vol(\sl(\O_2))\simeq1.003841$ \\
$vol(\sl(\O_3))\simeq0.169156$ \\
$vol(\sl(\O_7))\simeq0.888914$\\
$vol(\sl(\O_{11}))\simeq1.165895.$\\
\end{tabular}
\end{center}
Using the above volume formula, Borel proved in \cite{borel} that given $K>0$, there are only finitely many conjugacy classes of 
arithmetic lattices in $\sl(\C)$ such that volume of the associated 3-fold is $\leq K$. Work of Jorgensen and Thurston tell us 
that there is a hyperbolic 3-manifold of smallest volume and there are only finitely many non-isometric hyperbolic 3-manifolds 
attaining this minimal volume. It is conjectured that the so called {\em Weeks manifold}, of volume $\simeq 0.9427073$, attains this minimum. 
This manifold is in fact arithmetic given by a maximal order in a quaternion algebra, over the unique cubic number field $K$ of discriminant $-23$, which ramifies at the real places and the unique ideal of $K$ of norm $5$. Chinburg, Friedburg, Jones and Reid proved in \cite{chinburg} that the Weeks Manifold 
has the smallest volume among all arithmetic hyperbolic 3-manifolds. Recently Gabai, Meyerhoff and Milley (see \cite{gabai}) showed that this manifold is the smallest volume among all compact hyperbolic 3-manifolds. 

\section{Betti Numbers of Arithmetic Hyperbolic 3-Manifolds} 
A famous conjecture, known as the {\bf Virtual Betti Number Conjecture}, of Thurston says that every hyperbolic 3-manifold $M$ has a finite cover $M'$ whose first Betti number is positive, that is, $\dim H^1(M',\C) \geq 1$. A stronger version, known as the {\bf Virtual Infinite Betti Number Conjecture}, of this conjecture says that a hyperbolic 3-manifold has finite covers with arbitrarily large first Betti number. It is known that for arithmetic 3-folds (see results of Agol \cite{agol-betti}, Venkataramana \cite{venkataramana}, Cooper, Long and Reid \cite{cooper-long-reid}), the weak conjecture is equivalent to the stronger one. In the case of arithmetic 3-folds of congruence type, this was already proven by Borel \cite{borel-74} in 1974. 

\subsection{Bianchi groups}
Let $\Gamma$ be an arithmetic lattice in $\sl(\C)$ and $Y$ be its associated hyperbolic 3-fold. Assume that $Y$ is non-compact. As before, without loss of generality, let us assume that $\Gamma$ is a Bianchi group $\sl(\O_K)$. Let $X$ be the Borel-Serre compactification of $Y$ that we mentioned in the previous section. It is well-known that the rank of the image of the map $res:H^1(X,\C) \rightarrow H^1(\partial X,\C)$ is half of the rank of $H^1(\partial X, \C)$, giving that we have $\dim H^1_{Eis}(Y,\C) = h_K$. This already shows that the Virtual Betti Number Conjecture is trivially true for non-compact arithmetic hyperbolic 3-manifolds (except when $K=\Q(\sqrt{-d})$ with $d=1,3$, in which case, the cross-sections of the cusps are 2-orbifolds whose underlying manifolds are 2-spheres, thus giving $H^1_{Eis}(Y,\C)=0$).

It was observed by Bianchi \cite{bianchi} and Swan \cite{swan} in the several examples they computed that $\dim H^1_{cusp}(Y,\C) = 0$. Mennicke showed 
(see the end of \cite{swan}) that this was not true in general; he computed that for $K=\Q(\sqrt{-10})$, the dimension of $H^1_{cusp}(Y,\C)$ is $1$. Works of Zimmert \cite{zimmert}, Grunewald-Schwermer \cite{grunewald-schwermer-1}, Vogtmann \cite{vogtmann}, Rohlfs \cite{rohlfs} and Kr\"amer 
\cite{kraemer} showed that one has $\dim H^1_{cusp}(Y,\C) = 0$ only for  $K=\Q(\sqrt{-d})$ with 
$$d \in \{ 1, 2, 3, 5, 6, 7, 11, 15, 19, 23, 31, 39, 47, 71 \}.$$

Zimmert in his 1971 Bielefeld Diplomarbeit \cite{zimmert} showed that there is a set $Z(K)$, today called the {\em Zimmert Set}, depending on the imaginary quadratic field $K$ such that 
$$\dim H^1_{cusp}(Y, \C) \geq \#Z(K).$$
Grunewald and Schwermer \cite{grunewald-schwermer-1} modified the purely topological method of Zimmert and proved that 
$$\dim H^1_{cusp}(Y,\C) \rightarrow \infty$$
as the discriminant of $K$ grows, that is, $|\bigtriangleup_K| \rightarrow \infty$. Using the same method, they also proved in \cite{grunewald-schwermer-2}  
every finite index subgroup of a Bianchi group $\sl(\O)$ has finite index (torsion-free) subgroups with arbitrarily large cuspidal cohomology. Hence the class of non-compact arithmetic hyperbolic 3-folds also satisfies the Virtual Infinite Betti Number Conjecture. 

Now consider the case $Y$ is compact. Let us put the additional hypothesis that the defining field $K$ of $\Gamma$ is an imaginary quadratic  field. Then by the Jacquet-Langlands correspondence we have an injection 
$$H^1(\Gamma,\C) \hookrightarrow H^1_{cusp}(\Gamma',\C)$$
for some finite index subgroup $\Gamma'$ of the Bianchi group $\sl(\O_K)$. Now the above mentioned result of Grunewald and Schwermer guarantees that we can find a finite index subgroup $\Gamma''$ of $\Gamma$ such that there is an embedding 
$$H^1(\Gamma'',\C) \hookrightarrow H^1_{cusp}(\Gamma''',\C)$$
for some finite-index subgroup $\Gamma''' \subset \Gamma'$ such that the image is non-zero, proving that $Y$ satisfies the Virtual 
Betti Number Conjecture.

\subsection{General Results}
A general result of Clozel \cite{clozel} gives that if $\Gamma$ is an arithmetic lattice in $\sl(\C)$ with defining quaternion algebra $D/K$ such that for every finite place $v \in S(D)$, $K_v$ contains no quadratic extension of $\Q_p$, here $p|Nv$, then $\Gamma$ satisfies the Virtual 
Betti Number Conjecture. Note that this result covers all arithmetic $\Gamma$ such that $K$ is imaginary quadratic and all 
places $v \in S(D)$ are such that the ideal associated to $v$ is of residue degree one (in other words, $p$ splits in $K$, here $p | Nv$).

Another important result related to the virtual Betti numbers of arithmetic hyperbolic 3-folds is that of Labesse-Schwermer \cite{labesse-schwermer} which uses certain functoriality results of Langlands programme, such as {\em Base Change}, together with Jacquet-Langlands correspondence. More precisely, let $Y$ be a compact arithmetic hyperbolic 3-manifold with defining field $K$. Assume that there is proper subfield $F \subset K$ (necessarily totally real) such that there is tower 
$$ F=F_0 \subset F_1 \subset \hdots \subset F_m \subset F_{m+1}=K$$
of intermediate fields such that $F_{i+1}/ F_i$ a Galois extension of prime degree and cyclic Galois group or a cubic non-Galois 
extension. Then $Y$ satisfies the Virtual Betti Number Conjecture. Rajan \cite{rajan} has improved this result further to the following form: 
Assume that there is a proper subfield $F \subset K$ and a solvable extension $L/F$ such that $F \subset K \subset L$, then 
$Y$ satisfies the Virtual Betti Number Conjecture. 

Let us sketch the argument of Labesse and Schwermer in a simple case. Assume that $Y$ is a compact hyperbolic 3-manifold 
associated to an arithmetic lattice $\Gamma$ whose field of definition $K$ is a prime degree Galois extension over $\Q$ with 
Galois group $G=\langle \sigma \rangle$. By Jacquet-Langlands correspondence, there is a finite-index subgroup $\Gamma'$ of $\sl(\O_K)$ such that 
$$S_2(\Gamma) \hookrightarrow S_2(\Gamma'). $$
Now taking the base change lift of a suitable weight 2 cuspidal modular forms for some suitable congruence subgroup of $\sl(\Z)$, 
Labesse and Schwermer show that there is a finite index subgroup $\Gamma''$ of $\Gamma'$ (which is stabled under the Galois automorphism $\sigma$) such that $S_2(\Gamma'')^{\delta(D)\textrm{-new}}$ is non-trivial. Now Jacquet-Langlands 
correspondence together with the generalized Eichler-Shimura  
isomorphism implies that there is a finite index subgroup $\Gamma'''$ of $\Gamma$ such that $H^1(\Gamma''',\C) \not= 0$, 
proving that $Y$ satisfies the Virtual Betti Number Conjecture.  

\section{Rational Homology 3-Spheres}
In this section, we will discuss a result of N.Dunfield and F.Calegari \cite{calegari-dunfield} which answered a question of D.Cooper negatively. 
To present Cooper's question, let us introduce some terminology.  A 3-manifold $M$ is called a {\bf rational homology 3-sphere} if its first Betti number 
(that is, the dimension of $H^1(M,\Q)$) is zero. The {\bf injectivity radius} of $M$, denoted $injrad(M)$, is the radius of the largest ball that can be embedded around every point of $M$, equivalently, it is half of the length of the shortest closed geodesic in $M$. 

In the 1995 version of the Problems in Low Dimensional Topology edited by Kirby, Cooper asked the following: Does there exist 
a $K>0$ such that if $M$ is closed hyperbolic 3-manifold with $injrad(M)>K$, then $\dim H^1(M,\C) >0$ ? If the answer to this question is {\em yes}, the Virtual Betti Number Conjecture is true since one can always find finite covers with arbitrarily large injectivity radius. Assuming the Generalized Riemann Hypothesis and a standard conjecture in Langlands Programme, Dunfield and Calegari constructed a 
tower of {\em arithmetic} co-compact lattices $\Gamma_n$ such that $\dim H^1(\Gamma_n,\C)=0$ for all $n$ and 
$injrad(\Gamma_n) \rightarrow \infty$ as $n \rightarrow \infty$.   

Let us sketch the arguments. Start with the quaternion algebra $D=(\frac{-1,-3}{K})$ over $K=\Q(\sqrt{-2})$. Then $D$ only ramifies 
at the two finite places corresponding to the ideal $\p,\bar{\p}$ over $3$. Let $\O$ be a maximal order in $D$ and let $\Gamma_0=\slo(\O)$ be the co-compact arithmetic lattice in $\sl(\C)$.  Let $\Gamma_n$ be the congruence subgroup of $\Gamma_0$ 
of level $\p^n$. For large $n$, $\Gamma_n$ is torsion free and the injectivity radius of the associated manifolds can be proven, by an analysis of the traces of elements of $\Gamma_n$, to go to infinity as $n$ grows. 

The Jacquet-Langlands correspondence implies that 
$$H^1(\Gamma_n, \C) \hookrightarrow H^1_{cusp}(\Gamma_0(\bar{\p}) \cap \Gamma_1(\p^n),\C)$$
where $\Gamma_0(\bar{\p})$ and $\Gamma_1(\p^n)$ denote the standard congruence subgroups of the Bianchi group $\sl(\Z[\sqrt{-2}])$ which are defined respectively as the subgroup of elements which are of the form 
$(\begin{smallmatrix} * & * \\ 0 &* \end{smallmatrix})$ modulo $\bar{\p}$ and the subgroup of elements of the form 
$(\begin{smallmatrix} 1 & * \\ 0 &1 \end{smallmatrix})$ modulo $\p^n$. 

A standard conjecture in Langlands Programme claims the existence of a bijection between the classes in 
$H^1_{cusp}(\Gamma_0(\bar{\p}) \cap \Gamma_1(\p^n),\C)$ which are simultaneous eigenvectors under the action of the family 
$\mathbb{T}$ of Hecke operators and certain irreducible continuous representations of the absolute Galois group $G_K$ of 
$K$ into $\textrm{GL}_2(\overline{\Q}_3)$. As $G_K$ is compact, it can be shown that the images of these representations, 
after possibly conjugation, lie inside a lattice and thus one can reduce the image modulo $3$. Doing so, we get continuous 
representations of $G_K$ into $\textrm{GL}_2(\overline{\F}_3)$.  
In particular, these mod $3$ representions have finite image and hence their kernels 
cut out finite Galois extensions $L/K$ with special properties; namely the extension can ramify only over $3$ and its Galois group embeds into some $\textrm{GL}_2(\F_{3^a})$. Calegari and Dunfield prove the non-existence of such extensions $L/K$ under the Generalized Riemann Hypothesis and thus show, under the above standard conjecture from Langlands Programme, that $H^1_{cusp}(\Gamma_0(\bar{\p}) \cap \Gamma_1(\p^n),\C)=0$. From this we deduce that $H^1(\Gamma_n, \C)=0$ and 
thus the manifolds associated to $\Gamma_n$ are rational homology 3-spheres.   
 
\section{Fibrations of Arithmetic Hyperbolic 3-Manifolds}
In this section, we shall discuss an interesting result of N.Dunfield and D.Ramakrishnan \cite{dunfield-ramakrishnan} that relates the information coming from the action of Hecke operators on the cohomology of an arithmetic hyperbolic 3-manifold $M$ to the number of finite covers of $M$ which fiber over the circle.   

Let $S$ be a closed surface and $\phi$ be an self-homeomorphism of $S$. 
Then the {\em mapping torus with monodromy $\phi$} is the 3-manifold given by 
$$M_{\phi} := S \times [0,1] / \{ (x,0) \equiv (\phi(x),1) \}.$$
We get a fibration $M_{\phi} \rightarrow S^1$ over the circle. A deep result 
of Thurston \cite{thurston-arxiv} that $M_{\phi}$ is hyperbolic if and only if $\phi$ is of type {\em pseudo-Anasov}. The well-known 
{\bf Virtual Fibration Conjecture} of Thurston says that every hyperbolic 3-manifold has a finite cover which fibers over the circle. 
Recent result of I.Agol \cite{agol-fiber}, using work of Haglund and Wise \cite{haglund-wise}, shows that a large family of 
hyperbolic 3-folds satisfies this conjecture, including arithmetic hyperbolic 3-folds containing totally geodesic surfaces (this includes 3-folds associated to Bianchi groups). 

Prior to these developments, Dunfield and Ramakrishnan constructed an explicit example which satisfied a stronger form of the above conjecture. More precisely, they constructed an arithmetic compact hyperbolic manifold $M$ which has finite covers which fiber over the circle in arbitrarily large number of ways. 

If a 3-manifold $M$ fibers over the circle $S^1$, then the fibration $M \rightarrow S^1$ gives rise to a homomorphism at the level of fundamental groups 
$\pi_1(M) \rightarrow \Z$, which can be seen as an element of $H^1(M,\Z) \simeq H^1(\pi_1(M), \Z).$ Cohomology classes which arise from 
fibrations (over the circle) can be studied via the {\bf Thurston norm} on $H^1(M,\Z)$.  For a compact 
connected surface $S$, let $\chi^*(S)$ be the absolute value of the Euler characteristic $\chi(S)$ of $S$ if $\chi(S) \leq 0$, else 
let $\chi^*(S)$ be $0$. For a surface $S$ with connected components $S_i$, let $\chi^*(S)$ be the sum of 
$\chi^*(S_i)$. The Thurston norm, as introduced by Thurston in \cite{thurston1}, of $\phi \in H^1(M,\Z)$ is 
$$\norm{\phi} := \textrm{inf} \{ \chi^*(S) \mid S \ \textrm{is a properly embedded oriented surface that is dual to} \ \phi \}.$$
Let us explain the meaning of {\em dual}. Any oriented surface $S$ that is properly embedded in $M$ gives rise to a homology 
class $[S] \in H_2(M,\partial M, \Z)$. By Poincar'e duality, there is an isomorphism $PD : H_2(M, \partial M ,\Z) \rightarrow H^1(M,\Z)$ 
and classes mapping to each other under $PD$ are called duals of each other. The above definition extends continuously to 
$H^(M,\R)$ and gives in fact a {\em semi-norm}, that is, non-zero classes in $H^1(M,\R)$ may have zero norm. 

The unit ball $B_M$ in $H^1(M,\R)$ with respect to the Thurston norm is a finite polytope. A top dimensional face $F$ of $B$ which satisfies the 
following property is called a {\bf fibered face}:  $\omega \in H^1(M,\Z)$ arises from a fibration if and only if the ray from the origin 
through $\omega$ (here we see $H^1(M,\Z)$ as a lattice in $H^1(M,\R)$) intersects the interior of $F$.

The core idea of the paper of Dunfield and Ramakrishnan has two ingredients: a result of Fried on when two cohomology classes in $H^1(M,\Z)$ 
which arise from fibrations lie on distinct fibered faces, and an analysis of the behaviour of cohomology classes arising from fibrations under the 
action of Hecke operators. In a simplified form, their core lemma says the following. Let $M$ be a closed hyperbolic 3-manifold and $\omega \in H^1(M,\Z)$ a class arising from a fibration over the circle. Let $T_g$ be a Hecke operator associated to some $g \in Comm(\pi_1(M))$. If $T_g(\omega)=0$, then there is a finite cover $N$ of $M$ with two classes in $H^1(N,\Z)$ (which arise from $\omega$) which lie in distinct fibered faces of $B_N$.  

Now the goal is to find an arithmetic (so that we have an ample collection of Hecke operators) $M$ such that there is at least one class in 
$H^1(M,\Z)$ that arises from a fibration over the circle. To achieve this goal, one needs to be aware of the following result of Thurston \cite{thurston2}: 
let $S$ be a surface that is embedded in $M$ and assume that $S$ is totally geodesic. Then the dual in $H^1(M,\Z)$ of $[S] \in H^2(M,\Z)$ 
cannot arise from a fibration over the circle. 

Let us sketch how Dunfield and Ramakrishnan contruct their example. Let $f$ be the unique weight 2 cuspidal modular form for the 
arithmetic lattice $\Gamma_0(49) \subset \sl(\Z) \subset \sl(\R)$. Let $K=Q(\sqrt{-3})$ and $\p,\q$ be the two conjugate ideals in $\O_K$ over $7$. Let $F$ be the base change of $f$ to $K$. Then $F$ is weight 2 cuspidal modular form for the arithmetic lattice $\Gamma_0(\p^2\q^2) \subset 
\sl(\O_K) \subset \sl(\C)$. Let $D/K$ be the quaternion algebra with $S(D)=\{ \p,\q \}$. It can be shown that 
$F \in S_2(\Gamma_0(\p^2\q^2))^{\delta(D)\textrm{-new}}$ and thus there is an arithmetic cocompact lattice $\Gamma$ arising from some 
order in $D$ such that $F$ can be transferred to $S_2(\Gamma)$ by the Jacquet-Langlands correspondence. 

Let $M$ be the cocompact arithmetic hyperbolic 3-manifold associated to $\Gamma$. A difficult analysis shows that $H^1(M,\Q)$ 
has dimension 3 and the class corresponding to the image of $F$ under the Jacquet-Langlands map arises from a fibration of $M$ 
over the circle. Now we are ready to use Hecke operators to construct finite covers of $M$ whose Thurston norm unit balls have arbitrarily many fibered faces. It is well-known that there are infinitely many Hecke operators which act $0$ on $f \in S_2(\Gamma_0(49))$. It follows directly from properties of the base change operation that the same is true for $F \in S_2(\Gamma_0(\p^2\q^2))$. As Jacquet-Langlands map respects the Hecke action, we deduce that the same is true for the image of $F$ in $S_2(\Gamma) \simeq H^1(M,\R)$. Let $\mathcal{P}= \{ \a_i \}$ be the set of ideals $\a_i$ of $K$ for which the associated Hecke operators act as $0$ on $F$ (equivalently, the image of $F$ under the Jacquet-Langlands map). The core lemma 
of the paper now implies that congruence subgroups $\Gamma_n $ of $\Gamma$ of type $\Gamma_0$ and level $\a_1 \cdots \a_n$ have 
at least $2^n$ different fibered faces in the unit ball for the Thurston norm in their cohomology $H^1(\Gamma_n,\R)$ as desired.  

Finally let us give a motivation for the choice of the specific $f \in S_2(\Gamma_0(49))$ that Dunfield and Ramakrishnan worked with.  
Let $E/\Q$ is a quaternion algebra unramified at the unique real place of $\Q$ and ramified at $7$ and some prime $p$ which stays 
inert in $K$, then we have $E \otimes_\Q K \simeq D/K$ where $D/K$ is as above. In particular, for a suitable arithmetic lattice $\Sigma \subset \sl(\R)$, 
arising from $E/\Q$, the compact totally geodesic surface $S$ associated to $\Sigma$ embeds in the compact hyperbolic 3-manifolds 
$M$ above. If $f$ could be transfered to $E/\Q$, then the image of $F$ under the Jacquet-Langlands map 
in $S_2(\Gamma) \simeq H^1(M,\R)$ would be the dual of $[S] \in H_2(M,\Z)$ and thus could not be arising from a fibration of $M$ 
over the circle. Indeed, our specific $f$ above chosen by Dunfield and Ramakrishnan has the special property that it does not 
transfer to any quaternion algebra $E/\Q$ as above.

\end{document}